\documentclass{amsart}
\usepackage[T1]{fontenc}
\usepackage{amsmath,amssymb}
\usepackage{mathtools}
\usepackage{enumerate}
\usepackage{amsfonts} 
\usepackage{hyperref}
\usepackage{mathrsfs}
\usepackage{amsthm}
\usepackage{eufrak}
\usepackage{commath}
\usepackage{stackengine}
\usepackage{tikz}

\newtheorem{theorem}{Theorem}[section]

\newtheorem{lemma}[theorem]{Lemma}

\title[on the number of dot products in the plane]{an improved bound on the number of dot products determined by a finite point set in the plane}
\author{Michalis Kokkinos}
\address{Institute for Algebra, Johannes Kepler University, Linz, Austria}
\email{mike\_kokkinos@hotmail.com}
\date{}

\begin{document}

\begin{abstract}
We are interested to bound from below the number of distinct dot products determined by a finite set of points $P$ in the Euclidean plane. In this paper, we build on \cite{dot} to obtain the improved lower bound \[|\{p\cdot q : p,q\in P\}|\gtrsim |P|^{\frac2 3 + \frac{7}{1425}}.\]
\end{abstract}
\maketitle

\section{Introduction}
Let $P\subseteq\mathbb{R}^2$ be a finite set of points in the Euclidean plane. We define $\Lambda(P)$ to be the set of distinct dot products determined by the set $P$, that is\begin{equation*}
    \Lambda(P):=\{p\cdot q : p,q\in P\}.
\end{equation*}
The problem of obtaining a lower bound for the quantity $\Lambda(P)$ is inspired by the Erd\H{o}s distinct distances problem, which asks how many distinct Euclidean distances are determined by a finite set of $N$ points in the real plane. This question is solved (up to logarithmic factors) by L. Guth and N. Katz \cite{distance} who obtained the lower bound $N/\log N$. Similar bounds are believed to hold in the dot products problem we are investigating. In fact, we are not aware of any construction of a point set $P$ giving strictly less than order of $|P|$ distinct dot products.

Achieving a lower bound of order $|P|^\frac{2}{3}$ is a known simple application of the Szemer\'{e}di-Trotter Theorem  [Theorem 1, \cite{S-T}]. To see this, consider the set $T_c:=\{(p,q)\in P^2:p\cdot q=c\}$ of pairs of points from $P$ which determine a fixed dot product $c\in\mathbb{R}\!\setminus\!\{0\},$ and the set of lines $\mathcal{R}:=\{l_p:p\in P\}$, where $l_p$ is given by $(x,y)\in\mathbb{R}$ satisfying $xp_1+yp_2=c$, for some $p=(p_1,p_2)\in P.$ Then, by the Szemer\'{e}di-Trotter Theorem,\[|T_c|=\sum_{p\in P}|l_p\cap P|=|\{(p,l)\in P\times\mathcal{R}:p\in l\}|\ll|P|^\frac{4}{3},\]whence\[|P|^2\ll\!
\sum_{c\in\Lambda(P)\setminus\{0\}}|T_c|\ll|\Lambda(P)||P|^\frac{4}{3}\]rearranges to $|\Lambda(P)|\gg|P|^\frac{2}{3}.$

The first result \cite{dot} on this problem, due to B. Hanson, O. Roche-Newton, and S. Senger, pushes the exponent to a positive constant beyond the $2/3$ threshold bound, namely to $2/3+1/3057$. The methods used in \cite{dot}, and subsequently to our paper, are specific to the real setting of the problem. It is worth mentioning that the finite field and complex setting variants still remain at their threshold lower bound of $|P|^\frac{2}{3}$. In the first case, the threshold bound follows from a paper of M. Rudnev [Theorem 13, \cite{dotFp}] for sufficiently small sized point sets compared to the characteristic of the field, and in the latter case, from the Szemer\'{e}di-Trotter Theorem in the complex plane \cite{Toth} and a similar argument to the one demonstrated above. If we are instead interested in lower bounding the number of distinct wedge products rather than dot products determined by a finite point set, a constant beyond 2/3 in the exponent was proved both in the real and complex plane in \cite{bilinear} and in positive characteristic in \cite{wedgeFp}.

In this paper, we are focusing on the initial problem in the Euclidean plane and we squeeze the techniques and ideas used in \cite{dot} to obtain an improvement to the constant on the exponent.
As a result, we prove the following bound. 
\begin{theorem}\label{main}
    Let $P\subseteq\mathbb{R}^2$ be a finite set of points in the plane. Then\begin{equation*}
        |\Lambda(P)|\gtrsim|P|^{\frac{2}{3}+\frac{7}{1425}}.
    \end{equation*}
\end{theorem}

\section{Notation}
Throughout this paper we adopt the \textit{Vinogradov asymptotic notation.} Given variable quantities $X$ and $Y$, we write $X\ll Y$, or $Y\gg X$, to denote that there exists an absolute constant $c>0$ such that $X\le cY$. Furthermore, we write $X\gtrsim Y$, or $Y\lesssim X,$ to denote $X\gg Y/(\log_2(Y))^c$ for some absolute constant $c>0$. 
\section{Outline of the ideas}
Let $\mathcal{L}$ denote the set of lines through the origin covering $P$. In the forthcoming proof of Theorem \ref{main}, we restrict our attention to a certain configuration of points and lines in the plane, which yields the fewest distinct dot products among point sets $P$. In particular, according to Lemma \ref{pin}, we do not obtain more than order of $|P|^\frac{2}{3}$ distinct dot products, only if $|\mathcal{L}|$ is approximately equal to $|P|^\frac{1}{3}.$ Henceforth, we suppose that there exist about $|P|^\frac{1}{3}$ lines passing through the origin covering $P$, with each line containing about $|P|^\frac{2}{3}$ points of $P.$ As a further condition, we assume that each line contains a point of $P$ with $x$-coordinate $x_0$. To ensure all these assumptions, we need to pass to convenient subsets $P'\subseteq P$ of points, and $\mathcal{L''}\subseteq\mathcal{L}$ of lines. We achieve this without significant quantitative loss, as the sizes of these newly introduced subsets are very similar to the initial ones. 

Now we can utilise this framework to provide a simplified, albeit instructive, sketch of how we can lower bound $|\Lambda(P)|$. Denote the set containing all the slopes of lines from $\mathcal{L''}$ by $S$. Fix $s_0\in S$ and define $X$ to be the set of $x$-coordinates of points of $P'$ lying on the line with slope $s_0$. Then,\begin{equation*}
    \begin{split}
        |\Lambda(P)|&\ge|\{(x,s_0x)\cdot(x_0,sx_0):x\in X,s\in S\}|\\
        &=|\{xx_0+xx_0s_0s:x\in X,s\in S\}|\\
        &=|X(1+s_0S)|.
    \end{split}
\end{equation*}We want the size of this product set to grow at least as fast as $|P|^{\frac{2}{3}+c}$, for some constant $c>0$.

In reality, it is insufficient to consider only one fixed slope, so we rather fix two distinct elements of $S$, say $s$ and $s'$, and consider their two corresponding subsets of dot products. This is done to ensure the required growth of the quantity $|\Lambda(P)|$ in question, as we have no information about the size of $X(1+s_0S)$ growing sufficiently fast. We do have information, however, about the growth of a ratio of products of the shifts $(sS+1)$ and $(s'S+1)$, see Theorem \ref{super2}, and we are able to efficiently lower bound an intersection, $Z,$ of two subsets of $X$ after scaling, instead of working directly with the set $X$. In fact, achieving a greater lower bound for this intersection is a key step in improving the value of $c$ compared to previous papers, and it is relied on the observation that we can exploit the forthcoming Lemma \ref{tao} in this proof. The new trick we implement now is to directly relate the quantities $|Z(sS+1)|$ and $|Z(s'S+1)|$ with the expander (\ref{expander}), while suffering the least loss on the exponent as possible, by utilising carefully the Ruzsa Triangle and the Pl\"{u}nnecke-Ruzsa inequalities.

\begin{figure}[h]
\centering

\begin{tikzpicture}[scale=1.2]

\draw[->] (0,0) -- (5,0) node[right] {};
\draw[->] (0,0) -- (0,5) node[above] {};

\def\xzero{1.0}

\draw[dotted, thick] (\xzero,0) -- (\xzero,5);

\node[below] at (\xzero,0) {$x_0$};

\foreach \m in {0.5, 1.2, 2} {
  
  \pgfmathsetmacro{\xmax}{5/\m}
  
  \pgfmathsetmacro{\endpoint}{min(\xmax,5)}

  \draw[thick] (0,0) -- (\endpoint,{\m*\endpoint});
  
  \foreach \i in {0,1,2,3} {
      \pgfmathsetmacro{\x}{\xzero + \i*1.0}
      \ifdim \x pt < \endpoint pt
         \fill (\x,{\m*\x}) circle (2pt);
      \fi
  }
}
\end{tikzpicture}
\caption{Representation of $|P|^\frac{1}{3}$ lines through the origin, with each line containing $|P|^\frac{2}{3}$ points of $P$ including a point with $x$-coordinate $x_0.$ }
\end{figure}

We note that if one carries out the calculations, they will observe that due to the amount of variables in our expander (\ref{expander}), a quadratic growth or an even slower growth in Theorem \ref{super2} does not guarantee the push of the exponent of $|\Lambda(P)|$ beyond the $2/3$ threshold mark. In other words, some input from superquadratic expansion is crucial for this argument to work.

\section{Prerequisites}
The following two results play a crucial role in our analysis and constitute the core of the proof. The first one is a result concerning the growth of the ratio of products of certain shifts of a finite set of real numbers $X$. Showing that this expander grows faster than $|X|^{\frac{5}{2}}$ was the content of Section 2 in \cite{dot} and entails ideas such as the squeezing technique\footnote{illustrated in the proof of Theorem 2.1, in \cite{dot}} and similar elementary methods to \cite{iterated}. In this paper, we use an improvement of this bound obtained in \cite{super} which had reinforced these techniques with the Elekes-Szab\'{o} Theorem and computer algebra. We will refer to this result as the \textit{superquadratic expander}. Let $A^{(k)}:=\{a_1\dotsm a_k:a_1,\dotsc,a_k\in\ A\}$ denote the \textit{k-fold product set} of a real finite set $A$.
\begin{theorem}\label{super2}
    For any finite set $X\subseteq\mathbb{R}$, there exist $x,x'\in X$ such that \begin{equation}\label{expander}
        \Bigg|\frac{(xX+1)^{(2)}(x'X+1)^{(2)}}{(xX+1)^{(2)}(x'X+1)}\Bigg|\gtrsim|X|^\frac{31}{12}.
    \end{equation}
\end{theorem}The next result is concerned with pinned dot products in the plane. We omit its proof which is a direct application of the Szemer\'{e}di-Trotter Theorem and can be found in \cite{dot}. 
\begin{lemma}\label{pin}    Let $P$ be a finite set of points in $\mathbb{R}^2$. Let $\mathcal{L}$ be the set of all lines through the origin covering $P$. Then, there exists $p\in P$ such that\[|\{p\cdot q:q\in P\}|\gg|P|^{\frac{1}{2}}|\mathcal{L}|^{\frac{1}{2}}.\]
\end{lemma}

We make use of the standard Ruzsa Triangle inequality and a modified Pl\"{u}nnecke-Ruzsa inequality in the multiplicative setting. The former can be found in [Section 2.7,\cite{taovu}] and the latter as Lemma 1.5 in \cite{katzshen} or Corollary 2.4 in \cite{plunineq}.
\begin{lemma}\label{triangineq}
    Let $X,Y,Z$ be finite subsets in a commutative group. Then\begin{equation*}
        \bigg|\frac{Y}{Z}\bigg|\le\frac{|X Y||X Z|}{|X|}.
    \end{equation*}
\end{lemma}
\begin{lemma}\label{plunnecke}
    Let $X,B_1,\dotsc,B_k$ be finite subsets in a commutative group. Then, there exists $X'\subseteq X$, with $|X'|\ge|X|/2,$ such that \begin{equation*}
        |X'B_1\dotsm B_k|\le2^k\frac{|XB_1||XB_2|\dotsm|XB_k|}{|X|^{k-1}}.
    \end{equation*}
\end{lemma}
Furthermore, we apply the following inequality reflecting that pairs of sets with large multiplicative energy have large intersection after scaling. The additive analogue of the inequality can be found as Corollary 2.10 in \cite{taovu}.
\begin{lemma}\label{tao}
    Let $A, B$ be finite subsets in a commutative group. Then there exists $x\in A/B$ such that\begin{equation*}
        |A\cap xB|\ge \frac{E_*(A,B)}{|A||B|}\ge\frac{|A||B|}{|AB|},
    \end{equation*}where $E_*(A,B)$ denotes the \textit{multiplicative energy} of $A$ and $B$.
\end{lemma}

\section{Proof of Theorem \ref{main}}\label{proof}
Let $\mathcal{L}$ denote the set of lines passing through the origin and incident to $P$. We will show that \[|\Lambda(P)|\gtrsim|P|^{\frac{2}{3}+c},\]for $c=7/1425$. Furthermore, we assume \begin{equation}\label{assumption}
    |\mathcal{L}|\le|P|^{\frac{1}{3}+2c}
\end{equation}otherwise the result follows by Lemma \ref{pin}. 

It would be convenient to assume that each line in $\mathcal{L}$ passes through roughly the same number of points in $P$. We achieve this using dyadic decomposition at the cost of a logarithmic factor. Consider\begin{equation*}
    |P|=\sum_{l\in\mathcal{L}}|l\cap P|=\sum_{j=1}^{\lceil \log_2|P|\rceil}\underset{2^{j-1}\le|l\cap P|<2^j}{\sum_{l\in\mathcal{L}:}}|l\cap P|.
\end{equation*}
By the Pigeonhole Principle, there exists a $j_0$ for which the summand contributes at least the average to the sum. Write $M=2^{j_0-1}$ and let $\mathcal{L'}\subseteq\mathcal{L}$ be the set of lines satisfying $M\le|l\cap P|\le2M.$ Therefore,\begin{equation*}
    |P|\lesssim\frac{|P|}{\log|P|}\ll\sum_{l\in\mathcal{L'}}|l\cap P|\ll M|\mathcal{L'}|.
\end{equation*}
Importantly, observe that each line $l\in\mathcal{L'}$ contains at least\begin{equation*}
    M\gtrsim|P|^{\frac{2}{3}-2c}
\end{equation*}points of $P$, by (\ref{assumption}). We may use a rotation of the axes in order to assume that all the lines (up to a multiplicative constant) contained in $\mathcal{L'}$ have non-negative slope and the line lying on the $x$-axis is contained in our line set. Moreover, denote by $P'\subseteq P$ the subset of points covered by $\mathcal{L'}$. Note that $P'$ is a large subset of $P$, and in fact has the same size as $P$ up to a multiplicative constant and a logarithmic factor, as $|P'|\ge M|\mathcal{L'}|\gtrsim|P|$. 

The next step is to further restrict our attention to another large subset $\mathcal{L''}\subseteq\mathcal{L'}\subseteq\mathcal{L}$ containing the lines passing through points of $P'$ with some fixed $x$-coordinate. We do this because the set of pinned dot products determined by a point of $P'$ that lies on the $x$-axis is in bijection with the set of all vertical lines covering $P'$, which is in turn in bijection with the set $X$ defined to contain all $x$-coordinates of points of $P'$. We may therefore assume $|X|\le|P|^{\frac{2}{3}+c}$. For each $u\in X$, define\begin{equation*}
    \mathcal{L}(u):=\{l\in\mathcal{L'}:l\cap\{x=u\}\in P'\},
\end{equation*}and consider\begin{equation*}
    \sum_{u\in X}|\mathcal{L}(u)|=|P'|\gtrsim|P|.
\end{equation*}
Again, by the Pigeonhole Principle, there exists some $x_0\in X$ for which\begin{equation*}
    |\mathcal{L}(x_0)|\gtrsim\frac{|P|}{|X|}\ge|P|^{\frac{1}{3}-c}.
\end{equation*}Henceforth, define $\mathcal{L''}:=\mathcal{L}(x_0)$.

Now that we have passed to our convenient framework involving the subset $P'\subseteq P$ of points and $\mathcal{L''}\subseteq\mathcal{L}$ of lines, we need to introduce some additional notation. Namely, let us denote by $S$ the set of slopes of lines from the line set $\mathcal{L''}$ and apply Theorem \ref{super2} to find two elements $s,s'\in S$ such that we have \begin{equation}\label{superq}
    \Bigg|\frac{(sS+1)^{(2)}(s'S+1)^{(2)}}{(sS+1)^{(2)}(s'S+1)}\Bigg|\gtrsim|S|^\frac{31}{12}.
\end{equation}
Let $X(s)$ and $X(s')$ be the sets containing the $x$-coordinates of points of $P'$ lying on the lines of $\mathcal{L''}$ with slopes $s$ and $s'$ respectively, and for some $a\in\mathbb{R}\backslash\{0\}$, consider the intersection\begin{equation*}
    Z:=X(s)\cap aX(s').
\end{equation*}

We proceed by relating the quantities $|Z(sS+1)|$ and $|Z(s'S+1)|$ with $|\Lambda(P)|$. Consider the following two subsets of $\Lambda(P)$\begin{equation*}
    \begin{split}
        &\Lambda_1:=\{(x,sx)\cdot(x_0,s_1x_0):x\in Z, s_1\in S\},\\
        &\Lambda_2:=\{(x,s'x)\cdot(x_0,s_2x_0):x\in a^{-1}Z, s_2\in S\}.      
    \end{split}
\end{equation*}
In fact,\begin{equation*}
    \begin{split}
        &|\Lambda_1|=|\{(x,sx)\cdot(x_0,s_1x_0):x\in Z, s_1\in S\}|=|Z(sS+1)|,\\
        &|\Lambda_2|=|\{(x,s'x)\cdot(x_0,s_2x_0):x\in a^{-1}Z, s_2\in S\}|=|Z(s'S+1)|,   
    \end{split}
\end{equation*}whence\begin{equation}\label{relation}
    |Z(sS+1)|,|Z(s'S+1)|\le|\Lambda(P)|.
\end{equation}

The relations (\ref{superq}) and (\ref{relation}) may hint already that a strong lower bound on the intersection $Z$ will contribute to a stronger lower bound on the quantity $|\Lambda(P)|$ that we are interested in. To achieve this, consider the subset $\Lambda'\subseteq\Lambda(P)$ of all the dot products determined by the points of $P'$ lying on the lines of $\mathcal{L''}$ with slopes $s$ and $s'$. In particular\begin{equation*}
    \begin{split}
        \Lambda':\!&=\{(x,sx)\cdot(x',s'x'):x\in X(s), x'\in X(s')\}\\
        &=\{xx'(1+ss'):x\in X(s), x'\in X(s')\}\\
        &=X(s)X(s')(1+ss'),
    \end{split} 
\end{equation*}whence\begin{equation*}
    |X(s) X(s')|\le|P|^{\frac{2}{3}+c}.
\end{equation*}
Now, by Lemma \ref{tao}, there exists $a\in X(s)/X(s')$ for which\begin{equation*}
    |Z|=|X(s)\cap aX(s')|\ge\frac{|X(s)||X(s')|}{|X(s) X(s')|}\gg|P|^{\frac{2}{3}-5c},
\end{equation*}recalling each $X(s),X(s')$ contains roughly $M$ elements.

We put all the pieces together by implementing the aforementioned Ruzsa Triangle and Pl\"{u}nnecke-Ruzsa inequalities. Firstly, we make use of Lemma \ref{plunnecke} twice; there exists $Z_1\subseteq Z$, with $|Z_1|\ge|Z|/2,$ such that the inequality\begin{equation*}
    \big|Z_1(sS+1)^{(2)}(s'S+1)^{(2)}\big|\ll\frac{|Z(sS+1)|^2|Z(s'S+1)|^2}{|Z|^3}
\end{equation*}holds. Similarly, there exists $Z_2\subseteq Z_1$, with $|Z_2|\ge|Z_1|/2,$ for which we have\begin{equation*}
    \big|Z_2(sS+1)^{(2)}(s'S+1)\big|\ll\frac{|Z_1(sS+1)|^2|Z_1(s'S+1)|}{|Z_1|^2}.
\end{equation*}Note that $|Z_2|\gg|Z_1|\gg|Z|.$ Finally, we employ the superquadratic expander (\ref{superq}) and Lemma \ref{triangineq}, to deduce\begin{equation*}
    \begin{split}
        |S|^\frac{31}{12}&\lesssim\Bigg|\frac{(sS+1)^{(2)}(s'S+1)^{(2)}}{(sS+1)^{(2)}(s'S+1)}\Bigg|\\
        &\le\frac{\big|Z_2(sS+1)^{(2)}(s'S+1)^{(2)}\big|\big|Z_2(sS+1)^{(2)}(s'S+1)\big|}{|Z_2|}\\
        &\ll\frac{|Z(sS+1)|^4|Z(s'S+1)|^3}{|Z|^6}\\
        &\le\frac{|\Lambda(P)|^7}{|Z|^6}.
    \end{split}
\end{equation*}
Rearranging, and observing that $|S|=|\mathcal{L''}|\gg|P|^{\frac{1}{3}-c}$, yields\begin{equation*}
    \begin{split}
        |\Lambda(P)|&\gg|S|^\frac{31}{84}|Z|^{\frac{6}{7}}\\
        &\gg|P|^{\frac{25}{36}-\frac{391c}{84}}.
    \end{split}
\end{equation*}
We optimise by setting the exponent equal to $2/3+c$. This gives $c=7/1425$ and concludes the proof.

\hfill $\square$
\section{Acknowledgments}
I want to thank my supervisor Oliver Roche-Newton for his guidance and patience. His advices were proven to be invaluable for the outcome of this work. Special thanks to Audie Warren and Jakob F\"{u}hrer for many helpful discussions. I also thank the anonymous referees for their constructive comments and suggestions. This work was supported by the Austrian Science Fund PAT2559123.

\renewcommand{\bibname}{\Large References}
\bibliography{main}        
\bibliographystyle{plain}  

\end{document}